\documentclass[11pt, english]{article}
\usepackage[T2A]{fontenc}
\usepackage[cp1251]{inputenc}
\usepackage{babel}

\usepackage[centertags]{amsmath}
\usepackage{amsfonts}
\usepackage{amssymb}
\usepackage{amsthm}
\usepackage{newlfont}
\usepackage[active]{srcltx} 
\newtheorem{theorem}{\bf Theorem}
\newtheorem{lemma}{\bf Lemma}
\newtheorem{remark}{\bf Remark}

\newcommand{\mathP}{\mathbf P}
\newcommand{\mathE}{\mathbf E}

\newcommand{\einf}{\mathrm{ess\,inf\,}}

\begin{document}

\author{Andrei N. Frolov \footnote{ E-mail address: Andrei.Frolov@pobox.spbu.ru}
\\ St.~Petersburg State University }

\title{
Small deviations of iterated processes in space of trajectories
\footnotetext{ This work was partially supported by Federal
Program "Scientific and teaching staff of innovative Russia",
project 1.1-111-128-033.}
}
\maketitle

{\abstract{
We derive logarithmic asymptotics of probabilities of small
deviations for iterated processes in the space
of trajectories. We find conditions under which these asymptotics
coincide with those of processes generating iterated processes.
When these conditions fail the asymptotics are quite different.
}



}


\section{Introduction}

There are many investigations on asymptotic behaviour of probabilities of
small deviations for various classes of stochastic processes and sequences.
The most important studied classes are sums of independent random variables,
stochastic processes with independent increments and Gaussian processes.
One can find a detailed list of references in Lifshits \cite{L2}. 
We restrict our attention to iterated processes.

Let $\xi(t)$ and $\Lambda(t)$, $t\geqslant 0$, be independent stochastic processes,
defined on the same probability space. Assume that, with probability 1,
$\Lambda(t)$ has continuous trajectories, $\Lambda(t) \geqslant 0$ and
$\Lambda(0)=0$. The stochastic processes $\chi(t)= \xi(\Lambda(t))$,
$t\geqslant 0$, is called the iterated processes.

Let $y_t$ be a function such that $y_t \to \infty$ as $t \to \infty$.
Asymptotic behaviour of probabilities of small deviations
$$\mathP \left(\sup\limits_{0\leqslant u \leqslant t} |\chi(u)|
\leqslant y_t \right)$$
has been investigated by Frolov \cite{F3},
Martikainen, Frolov and Steinebach \cite{F4},
Aurzada and Lifshits \cite{AL5} and Frolov \cite{F6}
for  $\xi(t)$ and $\Lambda(t)$ from various classes of stochastic processes.
One can find results on logarithmic asymptotics
in these papers. We consider below the asymptotic behaviour of
small deviations in space of trajectories.

Assume that, with probability 1, trajectories of $\xi(t)$ are
right continuous functions and $\xi(0)=0$. Define a family
of processes as follows:

\begin{equation}\label{ksit}
\left\{ \xi_t(u) = \xi(t u),\; u \in [0,1],\; t\geqslant 0 \right\},
\end{equation}
where $y_t$ is a function such that
$y_t \rightarrow \infty$ as $t\rightarrow\infty$.
Put
$$Q_t =\mathP \left( \xi_t(\cdot) \in G y_t \right),$$
where $G\in \Re$, class $\Re$
of sets in the Skorohod space $D[0,1]$ is defined
in Mogul'skii \cite{M1} and $a G = \{a g: \, g \in G\}$ for $a \in \mathbb{R}$.

Mogul'skii \cite{M1} has studied the asymptotic behavior of
$\log Q_n$ as $n\rightarrow\infty$
for homogeneous process with independent increments $\xi(t)$ and
$\xi(t)=S_{[t]}$, where $S_{[t]}$ is a sum of $[t]$ independent,
identically distributed random variables.
Frolov \cite{F7} has considered a close problem for iterated compound Poisson
processes by studying of the asymptotics for probabilities
$$ R_t = \mathP \left( -g_1\left(\frac{\Lambda(u)}{\Lambda(t)}\right) y_t
\leqslant \chi(u) \leqslant g_2\left(\frac{\Lambda(u)}{\Lambda(t)}\right) y_t
\; \mbox{for all}\; u \in [0,t]\right).
$$
Note that using of random bounds is well motivated (see \cite{F3}, \cite{F7}
for details).
Probability $R_t$ turns 
to an analog of probability $Q_t$
for 
$$
\left\{ \chi_t(u) = \chi(t u),\; u \in [0,1],\; t\geqslant 0 \right\}
$$
in some special cases, only. If, for example, either $g_i(t)\equiv c_i$, or
$\Lambda(t)=\Lambda t^\beta$, where $\Lambda$ is a positive random variables
and $\beta>0$, then bounds $g_i\left(\frac{\Lambda(u)}{\Lambda(t)}\right) y_t$
in the definition of $R_t$ will be non-random. 
In general case, we will arrive at random set $G$.
This effect disappears if we consider family
\begin{equation}\label{etat}
\left\{ \eta_t(u) = \xi(\Lambda(t) u),\; u \in [0,1],\; t\geqslant 0
\right\}
\end{equation}
instead of $\{\chi_t(u)\}$ and put
$$ P_t = \mathP \left(\eta_t(\cdot) \in G y_t \right).
$$
In the sequel, we only consider
$G=\{g(t) \in D[0,1]: g(0)=0, -g_1(t) \leqslant g(t) \leqslant g_2(t)\}$,
where $g_i(t)$ are positive, continuous, non-decreasing functions such that
$g_i(0)>0$. Such $G$ are the most interesting and simple sets from $\Re$.
More general sets $G$ may be considered for some classes of processes.

In this paper,
we describe the asymptotic behavior of $\log P_t$ as $t \to \infty$.
We obtain generalizations of results in Frolov \cite{F3}, \cite{F7}, \cite{F6}.


\section{Results}

Let $\xi(t)$ and $\Lambda(t)$, $t\geqslant 0$, be independent stochastic processes,
defined on the same probability space. Assume that, with probability 1,
$\Lambda(t)$ has continuous trajectories, $\Lambda(t) \geqslant 0$,
$\Lambda(0)=0$, trajectories of $\xi(t)$ are
right continuous functions and $\xi(0)=0$. Define families
of processes $\left\{ \xi_t(u) \right\}$ and $\left\{ \eta_t(u) \right\}$
by relations (\ref{ksit}) and (\ref{etat}), correspondingly.

Let $g_i(t)$, $t \in [0,1]$, be positive, continuous, non-decreasing functions such that
$g_i(0)>0$, $i=1,2$. Put
$G=\{g(t) \in D[0,1]: g(0)=0, -g_1(t) \leqslant g(t) \leqslant g_2(t)\}$.
Let $\mathcal{G}$ be a set of subsets of $ D[0,1]$ consisting from
such $G$.
Denote $a G = \{ a g: \, g \in G\}$ for $a \in \mathbb{R}$.

Assume that
there exist positive functions $B(t)$ and $\zeta(t)$, $t>0$,
and positive functional H(G), $G\in\mathcal{G}$,
such that $B(t) \to \infty$ and $\zeta(t) \to 0$ as $t \to\infty$,
$\limsup\limits_{t \to\infty} B( c t)/B(t) < \infty$ for all
$ c >0$, and the relation
\begin{eqnarray}
&& \label{10}
\hspace*{-2\parindent}
\log \mathP \left(\xi_t(\cdot) \in G y_t \right)=
-   t H(G) \zeta(y_t) (1+o(1))
\quad \mbox{as} \quad t\to\infty,
\end{eqnarray}
holds for every positive function $y_t$ with
$y_t \to \infty$, $y_t=o(B(t))$ and $t \zeta(y_t)\to \infty$
as $t \to\infty$.

Put $\lambda_t=\einf\Lambda(t)$ and
$V_t(\lambda)=\mathP (\Lambda(t)<\lambda)$.

We start with a result in which the asymptotic of small deviations
of iterated processes is similar to that from (\ref{10}) and
$\lambda_t$ plays a role of $t$.

\medskip
\begin{theorem}\label{th1}
{\it
Assume that $\lambda_t \sim \tilde{\lambda}_t$ as $t\rightarrow\infty$,
where $\tilde{\lambda}_t$ is a continuous, strictly increasing
function. Suppose that $\lambda_t \to \infty$ as $t\rightarrow\infty$
and for every $c>1$
\begin{eqnarray}
&& \label{20}
\liminf\limits_{t\rightarrow\infty} V_t( c \lambda_t) >0.
\end{eqnarray}

Then
\begin{eqnarray}
&& \label{30}
\hspace*{-2\parindent}
\log \mathP \left(\eta_t(\cdot) \in G y_t \right) =
- \lambda_t H(G) \zeta(y_t) (1+o(1))
\quad \mbox{as} \quad t\to\infty
\end{eqnarray}
for every positive function $y_t$ with $y_t \to \infty$,
$y_t=o(B(\tilde{\lambda}_t))$ and $\lambda_t \zeta(y_t)\to \infty$
as $t \to\infty$.

}
\end{theorem}

Condition (\ref{20}) holds if, for example,
$\Lambda(t)$ have atoms of the same mass
in $\lambda_t$ for all $t$.

It may also happen that properties of $\Lambda(t)$
yields (\ref{20}).
The most simple and important example is
$\Lambda(t)=\Lambda f(t)$, where $\Lambda$ is a non-negative
random variable and $f(t)$ is a positive function such that
$f(t)\rightarrow\infty$ as $t\rightarrow\infty$.
The natural generalization of this example is
the following assumption.

Suppose that there exist a positive, increasing, continuous
function $f(t)$, $f(t)\rightarrow\infty$ as $t\rightarrow\infty$,
and a non-negative random variable $\Lambda$ such that
the distributions of $\Lambda(t)/f(t)$ converge weakly to
the distribution of $\Lambda$ as $t\rightarrow\infty$.
Note that $\Lambda$ may be degenerate.

Denote $\hat{\lambda}=\einf\Lambda$.

In our next result, we deal with
the case $\lambda_t=O(f(t))$ as $t\rightarrow\infty$.

\begin{theorem}\label{th2}
{\it
Assume that 
$\lambda_t/f(t) \to \hat{\lambda}$ as $t\rightarrow\infty$.

Then
\begin{eqnarray}
&& \label{40}
\hspace*{-2\parindent}
\log \mathP \left(\eta_t(\cdot) \in G y_t \right)=
- \hat{\lambda} f(t) H(G) \zeta(y_t) (1+o(1))
\quad \mbox{as} \quad t\to\infty
\end{eqnarray}
for every positive function $y_t$ with
$y_t \to \infty$, $y_t=o(B(f(t)))$
 and $f(t) \zeta(y_t)\to \infty$ as $t \to\infty$.

}
\end{theorem}

It is possible that $\lambda_t=o(f(t))$ as $t\rightarrow\infty$
in Theorem \ref{th2}. Moreover, it may happen that $\lambda_t=0$ for all $t>0$.
Hence (\ref{40}) turns to
$$\log \mathP \left(\eta_t(\cdot) \in G \right)
= o\left( f(t) \zeta(y_t)\right) \quad \mbox{as} \quad t\to\infty.
$$
Under additional assumptions, we have the following better result.

\begin{theorem}\label{th4}
{\it
Assume that $B(t)$ is a regularly varying at infinity function,
$\hat{\lambda}=0$ and $\lambda_t/f(t) \to \hat{\lambda}$ as $t\rightarrow\infty$.

Assume that for all $t>0$ the functions $F_t(\lambda)=V_t(\lambda f(t))$ and
$F(\lambda)=\mathP (\Lambda<\lambda)$ are continuous for
$\lambda \leqslant \lambda'$ and positive for
$\lambda\in(\lambda_t/f(t),\lambda']$ and $\lambda\in(0,\lambda']$,
correspondingly, where $\lambda' >0$.

Let $y_t$ be a positive function with
$y_t \to \infty$, $y_t=o(B(f(t)))$ and $f(t) \zeta(y_t)\to \infty$
as $t \to\infty$. Let $\varepsilon_t$ denote the solution of the equation
$$
\frac{-\log F_t(\varepsilon_t)}{\varepsilon_t} = f(t)  H(G) \zeta(y_t).
$$

Suppose that $\varepsilon_t f(t)$ is equivalent to a continuous, strictly
increasing function. Assume that
$y_t=o(B(\varepsilon_t f(t)))$, $\lambda_t=o(\varepsilon_t f(t))$
and for every $\tau>0$ the relation
\begin{eqnarray}
&& \label{60}
\hspace*{-2\parindent}
\log F_{t}(\tau \varepsilon_t) \sim
\log F_{t}(\varepsilon_t) 
\end{eqnarray}
holds as $t \to\infty$.

Then
\begin{eqnarray}
\log \mathP \left(\eta_t(\cdot) \in G y_t \right) =
\log F_{t}(\varepsilon_t) (1+o(1))
\label{70}
= - \varepsilon_t f(t)  H(G) \zeta(y_t) (1+o(1))
\quad \mbox{as} \quad t\to\infty.
\end{eqnarray}

Here $\varepsilon_t\rightarrow 0$ and $F_{t}(\varepsilon_t)\rightarrow 0$
as $t\rightarrow\infty$.

}
\end{theorem}

Theorem \ref{th4} implies that the asymptotic behaviour of
$\log \mathP \left(\eta_t(\cdot) \in G \right)$
may depend on properties of the distribution function
of $\Lambda$ at zero.

The following examples show that the righthand sides of
(\ref{40}) and (\ref{70}) may have different behaviours.

Assume that $\Lambda(t)=\Lambda f(t)$ and $\hat{\lambda}=0$.
If, for example, $F(\lambda) = \lambda^p$ for $\lambda\in[0,1]$,
$p>0$, then
$\log F_{t}(\varepsilon_t) \sim -p \log (f(t) \zeta(y_t))$
as $t\rightarrow\infty$.
If $F(\lambda) = (-\log \lambda)^{-p}$ for $\lambda\in(0,e^{-1}]$, $p>0$,
then
$\log F_{t}(\varepsilon_t) \sim -p \log\log (f(t) \zeta(y_t))$
as $t\rightarrow\infty$.

The following result yields that one can not omit the condition (\ref{60}).

\begin{theorem}\label{th5}
{\it Assume that all the conditions of Theorem \ref{th4} hold
besides the condition (\ref{60}). Assume that for every
$\tau>0$ the following relation holds
$\log F_{t}(\tau \varepsilon_t) \sim \tau^{p}
\log F_{t}(\varepsilon_t)$
as $t\rightarrow\infty$, where $p> 0$.

Then
\begin{eqnarray}
&& \label{80}\nonumber
\hspace*{-2\parindent}
\log \mathP \left(\eta_t(\cdot) \in G y_t \right)=
o( \varepsilon_t  f(t) \zeta(y_t) )
\quad \mbox{as} \quad t\to\infty.
\end{eqnarray}
}
\end{theorem}

Functions $B(t)$ and $\zeta(t)$ are usually related through
structures of considered processes in results on small deviations.
These relations are not used in proofs of Theorems \ref{th1}--\ref{th5}
and, therefore, 
we did not assume that they hold.
Nevertheless, we have the following results.

\begin{remark}\label{r1}
{\it
If there exist positive constants $d$ and $\gamma$ and a slowly
varying at infinity function $L(t)$ such that
$\zeta(t)=t^{-\gamma} L(t)$ and $t B^{-\gamma}(t) L(B(t)) \to d$
as $t\rightarrow\infty$, then one can omit the conditions
$f(t) \zeta(y_t)\to \infty$ and $y_t=o(B(\varepsilon_t f(t)))$
as $t\rightarrow\infty$  in Theorems \ref{th2}--\ref{th5} and
\ref{th4}--\ref{th5}, correspondingly.
}
\end{remark}

Turn to applications of Theorems \ref{th1}--\ref{th5}.
We consider relation (\ref{10}) 
as a result on asymptotic of small deviations in the space of trajectories.
Then sufficient conditions for (\ref{10}) may be taken
from known results or they may be derived
by applications of known technics.


Let $\{\eta_n\}$ be a sequence of independent, identically distributed random
variables. If $\mathE \, \eta_1$ exists, assume that $\mathE \, \eta_1 =0$.
Suppose that the distributions of
$(\eta_1+\eta_2+\cdots+\eta_n)/B_n$ converge weakly to a strictly stable
distribution $G_{\alpha}$, $\alpha\in (0,2]$, with
 $G_{\alpha}((-\infty,0)) \in (0,1)$,
where $\{ B_n \}$ is a sequence of positive constants.

Put $\xi(t)=\eta_1+\eta_2+\cdots+\eta_{[t]}$, $t\geqslant 0$.
By Theorem 1 in Mogul'skii \cite{M1}, the realtion
\begin{eqnarray*}
\log \mathP \left(\xi_t(\cdot) \in G y_t \right)=
- C H(G) \frac{t}{y_t^\alpha}  L(y_t) (1+o(1))
\quad \mbox{as} \quad t\to\infty
\end{eqnarray*}
holds for every positive function $y_t$ with
$y_t \to \infty$ and $y_t=o(B_{[t]})$ as $t \to\infty$, where
$$ H(G) = \int\limits_0^1 (g_2(t)+g_1(t))^{-\alpha} dt,
\quad L(x)= x^{\alpha-2}\,\mathE \, \eta_1^2 I\{ |\eta_1|< x\},
$$
$L(x)$ is a slowly varying at infinity function,
$C$ is an absolute constant, depending only on the distribution
$G_{\alpha}.$ If $\alpha=2$, then $C=\pi^2/2.$

Hence (\ref{10}) holds with $B(t)=B_{[t]}$ and $\zeta(t)= C t^{-\alpha} L(t)$
and the above results may be applied to
iterated processes generated by the sum process $\xi(t)$.
If, for example, $g_1(t)=g_2(t)\equiv 1$ we obtain Theorem 6 from Frolov \cite{F6}.

Applying Theorem 4 in Mogul'skii \cite{M1}, we arrive at
similar results for strictly stable processes $\xi(t)$ such that $\xi(1)$
has distribution $G_{\alpha}$.
Further applications of Theorems \ref{th1}--\ref{th5} may be derived
in the same way as it was done in Frolov \cite{F6} for $g_1(t)=g_2(t)\equiv 1$.


For some classes of stochastic processes, one may consider $G$ from wider sets
than $\mathcal{G}$. We permanently use in the proofs 
that probability $\mathP \left( \xi_\lambda(\cdot) \in G y_t \right)$
is a non-increasing function of $\lambda$.
We derive this monotonicity, supposing that $g_i(t)$ are non-decreasing.
The last condition on $g_i(t)$ may be omitted, if we assume that
this monotonicity holds at least for large $\lambda$.
For example, this assumption holds, if $\xi(t)$ is $H$-self-similar process
(i.e. finite dimensional distributions of $\xi(ct)$ and $c^{H}\xi(t)$
coincide for all $c>0$). Remember that fractional Brownian motion with Hurst
parameter $H$ is $H$-self-similar and strictly stable processes with
index $\alpha$ are $1/\alpha$-self-similar.


\section{Proofs}

We start with the following result.

\begin{lemma}\label{l2}
For every fixed $t$, probability
$\mathP \left( \xi_\lambda(\cdot) \in G y_t \right)$
is a non-increasing function of $\lambda$.
\end{lemma}

{\bf Proof.}
For $\lambda_2>\lambda_1>0$, we have
\begin{eqnarray}
\nonumber &&
\mathP \left( \xi_{\lambda_2}(\cdot) \in G y_t \right)
= \mathP \left( - g_1\left(\frac{v}{\lambda_2}\right) \leqslant
\frac{\xi( v)}{y_t} \leqslant g_2\left(\frac{v}{\lambda_2}\right)\;
\mbox{for all}\; v \in [0,\lambda_2] \right) \leqslant
\\ \nonumber &&
\mathP \left( - g_1\left(\frac{v}{\lambda_1}\right) \leqslant
\frac{\xi( v)}{y_t} \leqslant g_2\left(\frac{v}{\lambda_1}\right)\;
\mbox{for all}\; v \in [0,\lambda_2] \right) \leqslant
\\ \nonumber &&
\mathP \left( - g_1\left(\frac{v}{\lambda_1}\right) \leqslant
\frac{\xi( v)}{y_t} \leqslant g_2\left(\frac{v}{\lambda_1}\right)\;
\mbox{for all}\; v \in [0,\lambda_1] \right) =
\mathP \left( \xi_{\lambda_1}(\cdot) \in G y_t \right)
\end{eqnarray}
where we have used that $g_i(t)$ are non-decreasing. Note that
the last probability is 1 for $\lambda_1=0$.

\hfill $\Box$

We will also use the next result.

\begin{lemma}\label{l1}
Let $z_t$ be a positive function such that $z_t \rightarrow\infty$ and
$z_t\sim\tilde{z}_t$ as $t \rightarrow\infty$, where
$\tilde{z}_t$ is a continuous, strictly increasing, positive function.
If (\ref{10}) holds, then
\begin{eqnarray}
&& \label{31}
\hspace*{-2\parindent}
\log \mathP \left(\eta_{z_t}(\cdot) \in G y_t \right) =
- z_t H(G) \zeta(y_t) (1+o(1))
\quad \mbox{as} \quad t\to\infty
\end{eqnarray}
for every positive function $y_t$ such that $y_t \to \infty$ and
$y_t=o(B(\tilde{z}_t))$.
\end{lemma}

{\bf Proof.}
Let $y_t$ be a positive function such that $y_t \to \infty$ and
$y_t=o(B(\tilde{z}_t))$.

Let $\tilde{z}^{-1}_t$ be the inverse function to $z_t$.
Put $x_u = y_{\tilde{z}^{-1}_{u}}$ for $u>0$.
Hence $x_u = o(B(u))$ as $u \to \infty$.
By (\ref{10}), for every $\delta \in(0,1)$ there exists $U=U(\delta)$ such that
the inequality
$$
 e^{- \frac{1}{\delta} u H(G)  \zeta(x_u)} \leqslant
\mathP \left( \xi_u(\cdot) \in G x_u \right)
\leqslant e^{- \delta u H(G) \zeta(x_u)}
$$
holds for all $u>U$. Putting
$t=\tilde{z}^{-1}_{u}$ implies that
$$ e^{-\frac{1}{\delta}   \tilde{z}_t H(G) \zeta(y_t)} \leqslant
\mathP \left( \xi_{\tilde{z}_t}(\cdot) \in G y_t \right)
\leqslant e^{-\delta  \tilde{z}_t H(G) \zeta(y_t)}
$$
for all sufficiently large $t$.
It follows that
$$ - \frac{1}{\delta} \leqslant
\liminf \limits_{t \to \infty}
\frac{\log \mathP \left( \xi_{\tilde{z}_t}(\cdot) \in G y_t \right)}{  z_t H(G) \zeta(y_t)}
\leqslant
\limsup\limits_{t \to \infty}
\frac{\log \mathP \left( \xi_{\tilde{z}_t}(\cdot) \in G y_t \right)}{  z_t H(G) \zeta(y_t)}
\leqslant -\delta.$$
Passing to the limit as $\delta \rightarrow 1$
in the last relation, we get 
$$ \log \mathP \left( \xi_{\tilde{z}_t}(\cdot) \in G y_t \right) =
- z_t H(G) \zeta(y_t) (1+o(1)).
$$

Fix $c>0$. Using properties of $B(t)$, we conclude that
conditions $y_t=o(B(c \tilde{z}_t))$ and $y_t=o(B(\tilde{z}_t))$ are
equivalent. In the same way as before, we have that
$$ \log \mathP \left( \xi_{c \tilde{z}_t}(\cdot) \in G y_t \right) =
- c z_t H(G) \zeta(y_t) (1+o(1)).
$$
Assume now that  $c>1$. By Lemma \ref{l2},
$$ \mathP \left( \xi_{\frac{1}{c} \tilde{z}_t}(\cdot) \in G y_t \right)
\leqslant \mathP \left( \xi_{ z_t}(\cdot) \in G y_t \right)
\leqslant \mathP \left( \xi_{c \tilde{z}_t}(\cdot) \in G y_t \right)
$$
for all sufficiently large $t$. It follows that
$$  e^{-\frac{1}{c}   z_t H(G) \zeta(y_t)}
\leqslant \mathP \left( \xi_{ z_t}(\cdot) \in G y_t \right)
\leqslant  e^{-c   z_t H(G) \zeta(y_t)}
$$
for all sufficiently large $t$. The remainder of the proof is the same as that
for $\tilde{z}_t$ above.

\hfill $\Box$

{\bf Proof of Theorem \ref{th1}.} Let $y_t$ be a
positive function  with $y_t \to \infty$,
$y_t=o(B(\tilde{\lambda}_t))$ and $\lambda_t \zeta(y_t)\to \infty$
as $t \to\infty$.
Put
$ P_t = \mathP \left(\eta_t(\cdot) \in G y_t \right)$.
By the independence of $\xi(t)$ and $\Lambda(t)$, we have
\begin{eqnarray}
\nonumber
P_t =
\int\limits_{\lambda_t} ^{\infty} 
\mathP \left( \xi_\lambda(\cdot) \in G y_t \right)
\, dV_{t}(\lambda).
\end{eqnarray}

Take $\delta \in(0,1)$. By Lemma \ref{l2} and Lemma \ref{l1}
with $z_t = \lambda_t$, we get
$$ P_t \leqslant \mathP \left( \xi_{\lambda_t}(\cdot) \in G y_t \right)
\leqslant e^{-\delta  H(G) \lambda_t \zeta(y_t)}
$$
for all sufficiently large $t$. It yields that
$$ \limsup\limits_{t \to \infty} \frac{\log P_t}{ H(G) \lambda_t \zeta(y_t)}
\leqslant -\delta .$$
Passing to the limit as $\delta \rightarrow 1$,
we get the upper bound in (\ref{30}).

Now we turn to the lower bound.

Take $c>1$.
Using Lemma \ref{l2} and condition (\ref{20}), we have
\begin{eqnarray}
\nonumber
P_t \geqslant
\int\limits_{\lambda_t} ^{c \lambda_t}
\mathP \left( \xi_\lambda(\cdot) \in G y_t \right) \, dV_{t}(\lambda)
\geqslant
\mathP \left( \xi_{c \lambda_t}(\cdot) \in G y_t \right) \,
V_t(c\lambda_t)
\geqslant
L \,
\mathP \left( \xi_{c \lambda_t}(\cdot) \in G y_t \right)
\end{eqnarray}
for all sufficiently large $t$, where
$L=L(c)>0$.

Note that conditions $y_t=o(B(c \tilde{\lambda}_t))$ and
$y_t=o(B(\tilde{\lambda}_t))$ are equivalent.
Application of (\ref{31}) with $z_t=c \lambda_t$
yields that for every $\delta \in(0,1)$ the inequality
$$ \mathP \left( \xi_{c \lambda_t}(\cdot) \in G y_t \right)
\geqslant e^{-\delta c  \lambda_t H(G) \zeta(y_t)}
$$
holds for all sufficiently large $t$. Then
$$ \log P_t \geqslant \log L -\delta c  \lambda_t H(G) \zeta(y_t)$$
for all sufficiently large $t$. Taking into account that
$\lambda_t \zeta(y_t)\to \infty$
as $t \to\infty$, we arrive at the inequality
$$ \liminf\limits_{t \to \infty} \frac{\log P_t}{ \lambda_t H(G) \zeta(y_t)}
\geqslant - \delta c.
$$
Taking the limit as  $c \to 1$ and $\delta\to 1$,
we get the lower bound in (\ref{30}).

\hfill $\Box$

{\bf Proof of Theorem \ref{th2}.}
Assume first that $\hat{\lambda}>0$.
Put $\tilde{\lambda}_t=\hat{\lambda} f(t)$ and check (\ref{20}).

Take $c>1$. Chose $p \in (0,1)$ such that $p c>1$ and
$ p c \hat{\lambda}$ is the continuity point of
the distribution function $\Lambda$. We have
\begin{eqnarray}
\nonumber
\liminf\limits_{t\rightarrow\infty}
V_t(c \lambda_t)
= \liminf\limits_{t\rightarrow\infty}
\mathP \Big(\frac{\Lambda(t)}{f(t)} < c \frac{\lambda_t}{f(t)}\Big)
\geqslant
\liminf\limits_{t\rightarrow\infty}
\mathP \Big(\frac{\Lambda(t)}{f(t)} < p c \hat{\lambda}\Big)
= \mathP \big(\Lambda < p c \hat{\lambda}\big) > 0.
\end{eqnarray}

Theorem \ref{th1} yields that (\ref{40}) holds for
every positive function $y_t$ such that $y_t=o(B(\tilde{\lambda}_t))$
as $t \to\infty$. Taking into account that conditions
$y_t=o(B(\tilde{\lambda}_t))$ and $y_t=o(B(f(t)))$ are equivalent,
we finish the proof in the case $\hat{\lambda}>0$.

Turn to the case $\hat{\lambda}=0$.
Since $\log P_t \leqslant 0$,
we need only prove the lower bound.

Let $y_t$ be a positive function such that
$y_t \to \infty$,
$y_t=o(B(f(t)))$ and $f(t) \zeta(y_t)\to \infty$
as $t \to\infty$.

Take $\varepsilon>0$ such that $\varepsilon$ is a continuity point
of the distribution function of $\Lambda$. Since $\lambda_t=o(f(t))$
as $t \to \infty$ and $\mathP \left( \xi_\lambda(\cdot) \in G y_t \right)$
is a non-increasing function of $\lambda$, we have
in the same way as in the proof of Theorem \ref{th1} that
\begin{eqnarray*}
&& P_t \geqslant
\int\limits_{\lambda_t } ^{\varepsilon f(t) }
\mathP \left( \xi_\lambda(\cdot) \in G y_t \right) \, dV_t(\lambda)
\geqslant
\mathP \left( \xi_{\varepsilon f(t)}(\cdot) \in G y_t \right)
 \, V_t(\varepsilon f(t))
\end{eqnarray*}
for all sufficiently large $t$. Since
$V_t(\varepsilon f(t)) \to \mathP \big( \Lambda < \varepsilon) > 0$
as $t \to \infty$, then
\begin{eqnarray*}
&& P_t \geqslant C \,
\mathP \left( \xi_{\varepsilon f(t)}(\cdot) \in G y_t \right)
\end{eqnarray*}
for all sufficiently large $t$, where $C=C(\varepsilon)>0$.

Note that for every fixed $\varepsilon >0$
conditions $y_t=o(B(\varepsilon f(t)))$ and $y_t=o(B(f(t)))$
are equivalent. Making use of (\ref{31}) with $z_t=\varepsilon f(t)$, we get
$$ \mathP \left( \xi_{\varepsilon f(t)}(\cdot) \in G y_t \right)
\geqslant e^{- 2 \varepsilon f(t) \zeta(y_t)}
$$
for all sufficiently large $t$. It follows that
$$ \log P_t \geqslant \log C -2 \varepsilon f(t) \zeta(y_t)$$
for all sufficiently large $t$. Taking into account that
$f(t) \zeta(y_t)\to \infty$ as $t \to\infty$, we get
$$ \liminf\limits_{t \to\infty}
\frac{\log P_t}{f(t) \zeta(y_t)} \geqslant - 2 \varepsilon.
$$
Passing to the limit as $\varepsilon\rightarrow 0$,
we arrive at desired assertion.

\hfill $\Box$

{\bf Proof of Theorem \ref{th4}.}
Let $y_t$ be a function satisfying the conditions  of Theorem \ref{th4}.
Put $b_t=f(t) \zeta(y_t)$.
By assumptions, $b_t \rightarrow \infty$ as $t\rightarrow\infty$.
In the same way as in the proof of Theorem 3 on p. 172 in [3],
we have
$\varepsilon_t\rightarrow 0$ и $F_{t}(\varepsilon_t)\rightarrow 0$
as $t\rightarrow\infty$.

By the definition, $\varepsilon_t f(t) \rightarrow\infty$ as $t\rightarrow\infty$.
Note that for every fixed $c >0$ conditions
$y_t=o(B(\varepsilon_t f(t)))$ and $y_t=o(B(c\varepsilon_t f(t)))$
are equivalent.

Take $\delta\in(0,1)$. The inequality $\lambda_t \leqslant \varepsilon_t f(t)$
holds for all sufficiently large $t$. By Lemma \ref{l1} with
$z_t =\varepsilon_t f(t)$, we have
\begin{eqnarray}
\nonumber
&& P_t =
\int\limits_{\lambda_t} ^{\varepsilon_t f(t)}
\mathP \left( \xi_{\lambda}(\cdot) \in G y_t \right)
 \, dV_{t}(\lambda)+
\int\limits_{\varepsilon_t f(t)} ^{\infty}
\mathP \left( \xi_{\lambda}(\cdot) \in G y_t \right)
 \, dV_{t}(\lambda)
\\
&& \nonumber
\leqslant
F_t(\varepsilon_t)+
\mathP \left( \xi_{\varepsilon_t f(t)}(\cdot) \in G y_t \right)
= e^{-\varepsilon_t f(t)  H(G) \zeta(y_t)}+
\mathP \left( \xi_{\varepsilon_t f(t)}(\cdot) \in G y_t \right)
\leqslant 2 e^{-(1-\delta) \varepsilon_t f(t)  H(G) \zeta(y_t)}
\end{eqnarray}
for all sufficiently large $t$. This yields the upper bound in (\ref{70}).

Take $\tau>0$. The inequality $\lambda_t \leqslant \tau \varepsilon_t f(t)$
holds for all sufficiently large $t$. Applying Lemma \ref{l1} with
$z_t = \tau \varepsilon_t f(t)$, we have
\begin{eqnarray}
\nonumber
&& P_t \geqslant
\int\limits_{\lambda_t} ^{\tau \varepsilon_t f(t)}
\mathP \left( \xi_{\lambda}(\cdot) \in G y_t \right)
\, dV_{t}(\lambda)
\geqslant
\mathP \left( \xi_{\tau\varepsilon_t f(t)}(\cdot) \in G y_t \right)
\, F_t(\tau\varepsilon_t)
\\ \nonumber
&&
=e^{-\tau\varepsilon_t f(t)  H(G) \zeta(y_t)(1+o(1))} \, F_t(\tau\varepsilon_t)
=e^{-(1+\tau)\varepsilon_t f(t)  H(G) \zeta(y_t)(1+o(1))}.
\end{eqnarray}
This implies the lower bound in (\ref{70}).

\hfill $\Box$

{\bf Proof of Theorem \ref{th5}.} As in the proof of Theorem
\ref{th2} for $\hat{\lambda}=0$, we need only prove the lower bound.

Take $\tau>0$. In the same way as in the proof of Theorem \ref{th4}, we get
$$ P_t \geqslant
e^{-\tau\varepsilon_t f(t)  H(G) \zeta(y_t)(1+o(1))} \, F_t(\tau\varepsilon_t)
=
e^{-(\tau^{p}+\tau) \varepsilon_t f(t)  H(G) \zeta(y_t)(1+o(1))}
$$
as $t\rightarrow\infty$. The latter yields the lower bound.
\hfill $\Box$

One can find the proof of Remark \ref{r1} in Frolov \cite{F6}.


\end{document}